\documentclass[psamsfonts,english]{amsart}

    \usepackage[T1]{fontenc}
    \usepackage{babel,amsmath,amssymb,amsthm,url,mathtools,nicefrac}
        \mathtoolsset{showonlyrefs}
    \usepackage[mathcal]{eucal}

    \newtheorem{theorem}{Theorem}
    \newtheorem{lemma}[theorem]{Lemma}
    
    \newtheorem{algorithm}[theorem]{Algorithm}
    \theoremstyle{definition}
    \newtheorem{definition}[theorem]{Definition}
    \theoremstyle{remark}
    \newtheorem{remark}[theorem]{Remark}
    \newtheorem*{remark*}{Remark}

    \newcommand{\FF}{\mathbb{F}}
    \newcommand{\NN}{\mathbb{N}}
    \newcommand{\QQ}{\mathbb{Q}}
    \newcommand{\ZZ}{\mathbb{Z}}
    
    \newcommand{\CC}{\mathbb{C}}
    \newcommand{\PP}{\mathbb{P}}
    \newcommand{\jac}{\mathcal{J}_{C}}
    
    \newcommand{\frob}{\varphi}
    \newcommand{\basis}[1][B]{\mathcal{#1}}
    \newcommand{\pairing}{\varepsilon}
    \newcommand{\grp}[1]{\langle #1 \rangle}
    \newcommand{\zero}{\mathcal{O}}
    \newcommand{\pairingMatrix}[2][E]{\mathcal{#1}_{#2}}
    \newcommand{\C}{\mathcal{C}(\ell,q,k,\tau_k)}

    \DeclareMathOperator{\Div}{Div}
    
    \DeclareMathOperator{\Mat}{Mat}
    \DeclareMathOperator{\diag}{diag}
    \DeclareMathOperator{\rank}{rank}

    \newcommand{\nummerkontrol}[1]{}

\begin{document}

\title{Generators of Jacobians of Genus~Two~Curves}

\author[C.R. Ravnshøj]{Christian Robenhagen Ravnshøj}

\address{Department of Mathematical Sciences \\
University of Aarhus \\
Ny Munkegade \\
Building 1530 \\
DK-8000 Aarhus C}

\email{cr@imf.au.dk}

\thanks{Research supported in part by a PhD grant from CRYPTOMAThIC}




\begin{abstract}
We prove that in most cases relevant to cryptography, the Frobenius endomorphism on the Jacobian of a genus two curve
is represented by a diagonal matrix with respect to an appropriate basis of the subgroup of $\ell$-torsion points. From
this fact we get an explicit description of the Weil-pairing on the subgroup of $\ell$-torsion points. Finally, the
explicit description of the Weil-pairing provides us with an efficient, probabilistic algorithm to find generators of
the subgroup of $\ell$-torsion points on the Jacobian of a genus two curve.
\end{abstract}

\maketitle

\section{Introduction}

In \cite{koblitz87}, Koblitz described how to use elliptic curves to construct a public key cryptosystem. To get a more
general class of curves, and possibly larger group orders, Koblitz \cite{koblitz89} then proposed using Jacobians of
hyperelliptic curves. After Boneh and Franklin \cite{boneh-franklin} proposed an identity based cryptosystem by using
the Weil-pairing on an elliptic curve, pairings have been of great interest to cryptography~\cite{galbraith05}. The
next natural step was to consider pairings on Jacobians of hyperelliptic curves. Galbraith \emph{et
al}~\cite{galbraith07} survey the recent research on pairings on Jacobians of hyperelliptic curves.

Miller~\cite{miller} uses the Weil-pairing to determine generators of~$E(\FF_q)$, where~$E$ is an elliptic curve
defined over a finite field~$\FF_q$. Let $\jac$ be the Jacobian of a genus two curve defined over~$\FF_q$.
In~\cite{ravnshoj1}, the author describes an algorithm based on the Tate-pairing to determine generators of the
subgroup~$\jac(\FF_q)[m]$ of points of order~$m$ on the Jacobian, where~$m$ is a number dividing~$q-1$. The key
ingredient of the algorithm is a ``diagonalization'' of a set of randomly chosen points
$\{P_1,\dots,P_4,Q_1,\dots,Q_4\}$ on the Jacobian with respect to the (reduced) Tate-pairing~$\pairing$; i.e. a
modification of the set such that $\pairing(P_i,Q_j)\neq 1$ if and only if $i=j$. This procedure is based on solving
the discrete logarithm problem in~$\jac(\FF_q)[m]$. Contrary to the special case when $m$~divides~$q-1$, this is
infeasible in general. Hence, in general the algorithm in~\cite{ravnshoj1} does not apply.

In the present paper, we generalize the algorithm in~\cite{ravnshoj1} to subgroups of points of prime order~$\ell$,
where $\ell$ does not divide $q-1$. In order to do so, we must somehow alter the diagonalization step. We show and
exploit the fact that the $q$-power Frobenius endomorphism on~$\jac$ has a diagonal representation on~$\jac[\ell]$.
Hereby, computations of discrete logarithms are avoided, yielding the desired altering of the diagonalization step.

\subsection*{Setup}

Consider a genus two curve~$C$ defined over a finite field~$\FF_q$. Let $\ell$ be an odd prime number dividing the
number of $\FF_q$-rational points on the Jacobian~$\jac$, and with~$\ell$ dividing neither~$q$ nor~$q-1$. Assume that
the $\FF_q$-rational subgroup~$\jac(\FF_q)[\ell]$ of points on the Jacobian of order~$\ell$ is cyclic. Let~$k$ be the
multiplicative order of~$q$ modulo~$\ell$. Write the characteristic polynomial of the $q^k$-power Frobenius
endomorphism on~$\jac$ as
    $$P_k(X)=X^4+2\sigma_k X^3+(2q^k+\sigma_k^2-\tau_k)X^2+2\sigma_k q^kX+q^{2k},$$
where~$2\sigma_k,4\tau_k\in\ZZ$. Let $\omega_k\in\CC$ be a root of~$P_k(X)$. Finally, if $\ell$ divides~$4\tau_k$, we
assume that~$\ell$ is unramified in~$\QQ(\omega_k)$.

\begin{remark*}
Notice that in most cases relevant to cryptography, the considered genus two curve~$C$ fulfills these assumptions. Cf.
Remark~\ref{remark:ConditionsFulfilled} and~\ref{remark:AlmostAlwaysCyclic}.
\end{remark*}

\subsection*{The algorithm}

First of all, we notice that in the above setup, the $q$-power Frobenius endomorphism~$\frob$ on~$\jac$ can be
represented on~$\jac[\ell]$ by a diagonal matrix with respect to an appropriate basis~$\basis$ of~$\jac[\ell]$;
cf.~Theorem~\ref{theorem:FrobeniusDiagonalizable}. (In fact, to show this we do not need the $\FF_q$-rational
subgroup~$\jac(\FF_q)[\ell]$ of points on the Jacobian of order~$\ell$ to be cyclic.) From this observation it follows
that all non-degenerate, bilinear, anti-symmetric and Galois-invariant pairings on~$\jac[\ell]$ are given by the
matrices
    $$\pairingMatrix{a,b}=\begin{bmatrix}
            0 & a & 0 & 0 \\
            -a & 0 & 0 & 0 \\
            0 & 0 & 0 & b \\
            0 & 0 & -b & 0
            \end{bmatrix},\qquad a,b\in(\ZZ/\ell\ZZ)^\times
    $$
with respect to~$\basis$; cf.~Theorem~\ref{theorem:anti-symmetric-pairings}. By using this description of the pairing,
the desired algorithm is given as follows.

\setcounter{theorem}{16}

\begin{algorithm}
On input the considered curve~$C$, the numbers~$\ell$, $q$, $k$ and~$\tau_k$ and a number~$n\in\NN$, the following
algorithm outputs a generating set of~$\jac[\ell]$ or ``failure''.
\begin{enumerate}
 \item If $\ell$ does not divide $4\tau_k$, then do the following.
 \begin{enumerate}
  \item Choose points $\zero\neq x_1\in\jac(\FF_q)[\ell]$, $x_2\in\jac(\FF_{q^k})[\ell]\setminus\jac(\FF_q)[\ell]$ and~$x_3'\in
  U:=\jac[\ell]\setminus\jac(\FF_{q^k})[\ell]$; compute~$x_3=x_3'-\frob^k(x_3')$. If $\pairing(x_3,\frob(x_3))\neq 1$,
  then output~$\{x_1,x_2,x_3,\frob(x_3)\}$ and stop.
  \item Let $i=j=0$. While $i<n$ do the following
  \begin{enumerate}
   \item Choose a random point~$x_4\in U$.
   \item $i:=i+1$.
   \item If $\pairing(x_3,x_4)=1$, then $i:=i+1$. Else $i:=n$ and $j:=1$.
  \end{enumerate}
  \item If $j=0$ then output ``failure''. Else output~$\{x_1,x_2,x_3,x_4\}$.
 \end{enumerate}
 \item If $\ell$ divides $4\tau_k$, then do the following.
 \begin{enumerate}
  \item Choose a random point $\zero\neq x_1\in\jac(\FF_q)[\ell]$
  \item Let $i=j=0$. While $i<n$ do the following
  \begin{enumerate}
   \item Choose random points $y_3,y_4\in\jac[\ell]$; compute $x_\nu:=q(y_\nu-\frob(y_\nu))-\frob(y_\nu-\frob(y_\nu))$ for $\nu=3,4$.
   \item If $\pairing(x_3,x_4)=1$ then $i:=i+1$. Else $i:=n$ and $j:=1$.
  \end{enumerate}
  \item If $j=0$ then output ``failure'' and stop.
  \item Let $i=j=0$. While $i<n$ do the following
  \begin{enumerate}
   \item Choose a random point $x_2\in\jac[\ell]$.
   \item If $\pairing(x_1,x_2)=1$ then $i:=i+1$. Else $i:=n$ and $j:=1$.
  \end{enumerate}
  \item If $j=0$ then output ``failure''. Else output $\{x_1,x_2,x_3,x_4\}$ and stop.
 \end{enumerate}
\end{enumerate}
\end{algorithm}

\setcounter{theorem}{0}

Algorithm~\ref{algorithm:FindGenerators} finds generators of~$\jac[\ell]$ with probability at
least~$(1-\nicefrac{1}{\ell^n})^2$ and in expected running time~$O(\log\ell)$;
cf.~Theorem~\ref{theorem:FindGenerators}.

\begin{remark*}
To implement Algorithm~\ref{algorithm:FindGenerators}, we need to find a \emph{$q^k$-Weil number}
(cf.~Definition~\ref{definition:WeilNumber}). On Jacobians generated by the \emph{complex multiplication method}
\cite{weng03,gaudry,eisen-lauter}, we know the Weil numbers in advance. Hence, Algorithm~\ref{algorithm:FindGenerators}
is particularly well suited for such Jacobians.
\end{remark*}

\subsection*{Assumption}

In this paper, a \emph{curve} is an irreducible nonsingular projective variety of dimension one.

\section{Genus two curves}\label{sec:HyperellipticCurves}

A hyperelliptic curve is a projective curve $C\subseteq\PP^n$ of genus at least two with a separable, degree two
morphism~$\phi:C\to\PP^1$. It is well known, that any genus two curve is hyperelliptic. Throughout this paper, let~$C$
be a curve of genus two defined over a finite field~$\FF_q$ of characteristic~$p$. By the Riemann-Roch Theorem there
exists a birational map \mbox{$\psi:C\to\PP^2$}, mapping~$C$ to a curve given by an equation of the form
    $$y^2+g(x)y=h(x),$$
where $g,h\in\FF_q[x]$ are of degree $\deg(g)\leq 3$ and $\deg(h)\leq 6$; cf.~\cite[chapter~1]{cassels}.

The set of principal divisors $\mathcal{P}(C)$ on~$C$ constitutes a subgroup of the degree zero divisors~$\Div_0(C)$.
The Jacobian~$\jac$ of~$C$ is defined as the quotient
    $$\jac=\Div_0(C)/\mathcal{P}(C).$$
The Jacobian is an abelian group. We write the group law additively, and denote the zero element of the Jacobian
by~$\zero$.

{\samepage Let $\ell\neq p$ be a prime number. The $\ell^n$-torsion subgroup~$\jac[\ell^n]\subseteq\jac$ of points of
order dividing $\ell^n$ is a $\ZZ/\ell^n\ZZ$-module of rank four, i.e.
    \begin{equation*}\label{eq:J-struktur}
    \jac[\ell^n]\simeq\ZZ/\ell^n\ZZ\times\ZZ/\ell^n\ZZ\times\ZZ/\ell^n\ZZ\times\ZZ/\ell^n\ZZ;
    \end{equation*}
cf.~\cite[Theorem~6, p.~109]{lang59}.}

The multiplicative order~$k$ of~$q$ modulo~$\ell$ plays an important role in cryptography, since the (reduced)
Tate-pairing is non-degenerate over~$\FF_{q^k}$; cf.~\cite{hess}.

\begin{definition}[Embedding degree]
Consider a prime number $\ell\neq p$ dividing the number of $\FF_q$-rational points on the Jacobian~$\jac$. The
embedding degree of~$\jac(\FF_q)$ with respect to $\ell$ is the least number~$k$, such that~$q^k\equiv 1\pmod{\ell}$.
\end{definition}

\section{The Frobenius endomorphism}

Since~$C$ is defined over~$\FF_q$, the mapping $(x,y)\mapsto (x^q,y^q)$ is a morphism on~$C$. This morphism induces the
$q$-power Frobenius endo\-morphism $\frob$ on the Jacobian~$\jac$. Let $P(X)$ be the characteristic polynomial
of~$\frob$; cf.~\cite[pp.~109--110]{lang59}. $P(X)$ is called the \emph{Weil polynomial} of~$\jac$, and
    $$|\jac(\FF_q)|=P(1)$$
by the definition of~$P(X)$ (see~\cite[pp.~109--110]{lang59}); i.e. the number of~$\FF_q$-rational points on the
Jacobian is~$P(1)$.

\begin{definition}[Weil number]\label{definition:WeilNumber}
Let notation be as above. Let $P_k(X)$ be the characteristic polynomial of the $q^m$-power Frobenius
endomorphism~$\frob_m$ on~$\jac$. A~complex number~$\omega_m\in\CC$ with $P_m(\omega_m)=0$ is called a
\emph{$q^m$-Weil~number} of~$\jac$.
\end{definition}

\begin{remark}\label{rem:P_reducible}
Note that~$\jac$ has four $q^m$-Weil~numbers. If $P_1(X)=\prod_i(X-\omega_i)$, then $P_m(X)=\prod_i(X-\omega_i^m)$.
Hence, if~$\omega$ is a $q$-Weil number of~$\jac$, then~$\omega^m$ is a $q^m$-Weil number of~$\jac$.
\end{remark}

\section{Non-cyclic subgroups}\label{sec:NonCyclicSubgroups}

Consider a genus two curve~$C$ defined over a finite field~$\FF_q$. Let~$P_m(X)$ be the characteristic polynomial of
the $q^m$-power Frobenius endomorphism~$\frob_m$ on the Jacobian~$\jac$. $P_m(X)$ is of the form
$P_m(X)=X^4+sX^3+tX^2+sq^mX+q^{2m}$, where~$s,t\in\ZZ$. Let $\sigma=\frac{s}{2}$ and~$\tau=2q^m+\sigma^2-t$. Then
    $$P_m(X)=X^4+2\sigma X^3+(2q^m+\sigma^2-\tau)X^2+2\sigma q^mX+q^{2m},$$
and $2\sigma,4\tau\in\ZZ$. In \cite{ravnshoj2}, the author proves the following Theorem~\ref{theorem:rank}
and~\ref{theorem:rank:supplement}.

\begin{theorem}\label{theorem:rank}
Consider a genus two curve~$C$ defined over a finite field~$\FF_q$. Write the characteristic polynomial of the
$q^m$-power Frobenius endomorphism on the Jacobian~$\jac$ as
    $P_m(X)=X^4+2\sigma X^3+(2q^m+\sigma^2-\tau)X^2+2\sigma q^mX+q^{2m},$
where $2\sigma,4\tau\in\ZZ$. Let~$\ell$ be an odd prime number dividing the number of $\FF_q$-rational points on
$\jac$, and with $\ell\nmid q$ and~$\ell\nmid q-1$. If~$\ell\nmid 4\tau$, then
    \begin{enumerate}
    \item~$\jac(\FF_{q^m})[\ell]$ is of rank at most two as a $\ZZ/\ell\ZZ$-module, and
    \item~$\jac(\FF_{q^m})[\ell]$ is bicyclic if and only if $\ell$ divides $q^m-1$.
    \end{enumerate}
\end{theorem}

\begin{theorem}\label{theorem:rank:supplement}
Let notation be as in Theorem~\ref{theorem:rank}. Furthermore, let $\omega_m$ be a $q^m$-Weil number of~$\jac$, and
assume that~$\ell$ is unramified in~$\QQ(\omega_m)$. Now assume that~\mbox{$\ell\mid 4\tau$}. Then the following holds.
 \begin{enumerate}
 \item If $\omega_m\in\ZZ$, then $\ell\mid q^m-1$ and~$\jac[\ell]\subseteq\jac(\FF_{q^m})$.
 \item If $\omega_m\notin\ZZ$, then $\ell\nmid q^m-1$, $\jac(\FF_{q^m})[\ell]\simeq(\ZZ/\ell\ZZ)^2$ and
      ~$\jac[\ell]\subseteq\jac(\FF_{q^{mk}})$ if and only if $\ell\mid q^{mk}-1$.
 \end{enumerate}
\end{theorem}

Inspired by Theorem~\ref{theorem:rank} and~\ref{theorem:rank:supplement} we introduce the following notation.

\begin{definition}
Consider a curve~$C$ with Jacobian~$\jac$. We call~$C$ a~$\C$-curve, and write $C\in\C$, if the following holds.
    \begin{enumerate}
        \item $C$ is of genus two and defined over the finite field~$\FF_q$.
        \item $\ell$ is an odd prime number dividing the number of $\FF_q$-rational points on~$\jac$,
              $\ell$~divides neither~$q$ nor~$q-1$, and $\jac(\FF_q)$ is of embedding degree~$k$ with respect to~$\ell$.
        \item The characteristic polynomial of the
              $q^k$-power Frobenius endomorphism on~$\jac$ is given by
                $P_k(X)=X^4+2\sigma_k X^3+(2q^k+\sigma_k^2-\tau_k)X^2+2\sigma_k q^kX+q^{2k},$
              where~$2\sigma_k,4\tau_k\in\ZZ$.
        \item Let $\omega_k$ be a $q^k$-Weil number of~$\jac$.
              If $\ell$ divides~$4\tau_k$, then~$\ell$ is unramified in~$\QQ(\omega_k)$.
    \end{enumerate}
\end{definition}

\begin{remark}\label{remark:ConditionsFulfilled}
Since~$\ell$ is ramified in~$\QQ(\omega_k)$ if and only if~$\ell$ divides the discriminant of~$\QQ(\omega_k)$,
$\ell$~is unramified in~$\QQ(\omega_k)$ with probability approximately $1-\nicefrac{1}{\ell}$. Hence, in most cases
relevant to cryptography a genus two curve~$C$ is a~$\C$-curve.
\end{remark}

\section{Matrix representation of the Frobenius endomorphism}

An endomorphism $\psi:\jac\to\jac$ induces a linear map $\bar{\psi}:\jac[\ell]\to\jac[\ell]$ by restriction. Hence,
$\psi$ is represented by a matrix $M\in\Mat_4(\ZZ/\ell\ZZ)$ on~$\jac[\ell]$. If~$\psi$ can be represented on
$\jac[\ell]$ by a diagonal matrix with respect to an appropriate basis of~$\jac[\ell]$, then we say that~$\psi$ is
\emph{diagonalizable} or has a \emph{diagonal representation} on~$\jac[\ell]$.

Let $f\in\ZZ[X]$ be the characteristic polynomial of~$\psi$ (see~\cite[pp.~109--110]{lang59}), and let
$\bar{f}\in(\ZZ/\ell\ZZ)[X]$ be the characteristic polynomial of~$\bar{\psi}$. Then~$f$ is a monic polynomial of degree
four, and by \cite[Theorem~3, p.~186]{lang59},
    \begin{equation*}\label{eq:KarPolKongruens}
    f(X)\equiv\bar{f}(X)\pmod{\ell}.
    \end{equation*}

We wish to show that in most cases, the~$q$-power Frobenius endomorphism~$\frob$ is diagonalizable on~$\jac[\ell]$. To
do this, we need to describe the matrix representation in the case when~$\frob$ is not diagonalizable on~$\jac[\ell]$.

\begin{lemma}\label{lemma:FrobMatrix}\index{Frobenius endomorphism!matrix representation}
Consider a curve $C\in\C$. Let~$\frob$ be the $q$-power Frobenius endomorphism on the Jacobian~$\jac$. If~$\frob$ is
not diagonalizable on~$\jac[\ell]$, then~$\frob$ is represented on~$\jac[\ell]$ by a matrix of the form
    \begin{equation}\label{eq:FrobMatrix}
    M=\begin{bmatrix}
        1 & 0 & 0 & 0 \\
        0 & q & 0 & 0 \\
        0 & 0 & 0 & -q \\
        0 & 0 & 1 & c
        \end{bmatrix}
    \end{equation}
with respect to an appropriate basis of~$\jac[\ell]$.
\end{lemma}

\begin{proof}
Let $\bar{P}_k\in(\ZZ/\ell\ZZ)[X]$ be the characteristic polynomial of the restriction of the $q^k$-power Frobenius
endomorphism~$\frob_k$ to~$\jac[\ell]$. Since~$\ell$ divides the number of $\FF_q$-rational points on~$\jac$, $1$~is a
root of~$\bar{P}_k$. Assume that~$1$ is an root of~$\bar{P}_k$ with multiplicity~$\nu$. Then
    $$\bar{P}_k(X)=(X-1)^\nu\bar{Q}_k(X),$$
where $\bar{Q}_k\in(\ZZ/\ell\ZZ)[X]$ is a polynomial of degree~$4-\nu$, and~$\bar{Q}_k(1)\neq 0$. Since the roots
of~$\bar{P}_k$ occur in pairs~$(\alpha,1/\alpha)$, $\nu$~is an even number. Let $U_k=\ker(\frob_k-1)^\nu$
and~$W_k=\ker(\bar{Q}_k(\frob_k))$. Then~$U_k$ and~$W_k$ are $\frob_k$-invariant submodules of the
$\ZZ/\ell\ZZ$-module~$\jac[\ell]$, $\rank_{\ZZ/\ell\ZZ}(U_k)=\nu$, and~$\jac[\ell]\simeq U_k\oplus W_k$.

Assume at first that~$\ell$ does not divide~$4\tau_k$. Then~$\jac(\FF_q)[\ell]$ is cyclic and~$\jac(\FF_{q^k})[\ell]$
bicyclic; cf.~Theorem~\ref{theorem:rank}. By~\cite[Theorem~3.1]{rubin-silverberg02}, $\nu=2$. Choose
points~$x_1,x_2\in\jac[\ell]$, such that~$\frob(x_1)=x_1$ and $\frob(x_2)=qx_2$. Then~$\{x_1,x_2\}$ is a basis
of~$\jac(\FF_{q^k})[\ell]$. Now, let~$\{x_3,x_4\}$ be a basis of~$W_k$, and consider the
basis~$\basis=\{x_1,x_2,x_3,x_4\}$ of~$\jac[\ell]$. If~$x_3$ and~$x_4$ are eigenvectors of~$\frob_k$, then~$\frob_k$ is
represented by a diagonal matrix on~$\jac[\ell]$ with respect to~$\basis$. Assume~$x_3$ is not an eigenvector
of~$\frob_k$. Then~$\basis[B']=\{x_1,x_2,x_3,\frob_k(x_3)\}$ is a basis of~$\jac[\ell]$, and~$\frob_k$ is represented
by a matrix of the form~\eqref{eq:FrobMatrix}.

Now, assume~$\ell$ divides~$4\tau_k$. Since~$\ell$ divides~$q^k-1$, it follows
that~$\jac[\ell]\subseteq\jac(\FF_{q^k})$; cf.~Theorem~\ref{theorem:rank:supplement}. Let~$\bar{P}\in(\ZZ/\ell\ZZ)[X]$
be the characteristic polynomial of the restriction of~$\frob$ to~$\jac[\ell]$. Since~$\ell$ divides the number of
$\FF_q$-rational points on~$\jac$, $1$~is a root of~$\bar{P}$. Assume that~$1$ is an root of~$\bar{P}$ with
multiplicity~$\nu$. Since the roots of~$\bar{P}$ occur in pairs~$(\alpha,q/\alpha)$, it follows that
    $$\bar{P}(X)=(X-1)^\nu(X-q)^\nu\bar{Q}(X),$$
where $\bar{Q}\in(\ZZ/\ell\ZZ)[X]$ is a polynomial of degree $4-2\nu$, $\bar{Q}(1)\neq 0$ and~$\bar{Q}(q)\neq 0$. Let
$U=\ker(\frob-1)^\nu$, $V=\ker(\frob-q)^\nu$ and~$W=\ker(\bar{Q}(\frob))$. Then~$U$, $V$ and~$W$ are $\frob$-invariant
submodules of the $\ZZ/\ell\ZZ$-module~$\jac[\ell]$, $\rank_{\ZZ/\ell\ZZ}(U)=\rank_{\ZZ/\ell\ZZ}(V)=\nu$, and
$\jac[\ell]\simeq U\oplus V\oplus W$. If~$\nu=1$, then it follows as above that~$\frob$ is either diagonalizable
on~$\jac[\ell]$ or represented by a matrix of the form~\eqref{eq:FrobMatrix} with respect to some basis
of~$\jac[\ell]$. Hence, we may assume that~$\nu=2$. Now choose~$x_1\in U$, such that~$\frob(x_1)=x_1$, and expand this
to a basis $(x_1,x_2)$ of~$U$. Similarly, choose a basis $(x_3,x_4)$ of~$V$ with~$\frob(x_3)=qx_3$. With respect to the
basis~$\basis=\{x_1,x_2,x_3,x_4\}$, $\frob$~is represented by a matrix of the form
    $$
    M=\begin{bmatrix}
      1 & \alpha & 0 & 0 \\
      0 & 1 & 0 & 0 \\
      0 & 0 & q & \beta \\
      0 & 0 & 0 & q
      \end{bmatrix}.
    $$
Notice that
    $$
    M^k=\begin{bmatrix}
      1 & k\alpha & 0 & 0 \\
      0 & 1 & 0 & 0 \\
      0 & 0 & 1 & k q^{k-1}\beta \\
      0 & 0 & 0 & 1
      \end{bmatrix}.
    $$
Since~$\jac[\ell]\subseteq\jac(\FF_{q^k})$, we know that $\frob^k=\frob_k$ is the identity on~$\jac[\ell]$. Hence,
$M^k=I$. So~$\alpha\equiv\beta\equiv 0\pmod{\ell}$, i.e.~$\frob$ is represented by a diagonal matrix with respect
to~$\basis$.
\end{proof}

The next step is to determine when the Weil polynomial splits modulo~$\ell$.

\begin{lemma}\label{lemma:quadraticResidue}
Consider a curve $C\in\C$. Let~$\frob$ be the $q$-power Frobenius endomorphism on the Jacobian~$\jac$. Assume
that~$\frob$ is not diagonalizable on~$\jac[\ell]$, and let~$\frob$ be represented on~$\jac[\ell]$ by the matrix
    \begin{equation}
    M=\begin{bmatrix}
        1 & 0 & 0 & 0 \\
        0 & q & 0 & 0 \\
        0 & 0 & 0 & -q \\
        0 & 0 & 1 & c
        \end{bmatrix}
    \end{equation}
with respect to an appropriate basis of~$\jac[\ell]$. Let $P_n(X)$ be the characteristic polynomial of the $q^n$-power
Frobenius endomorphism on~$\jac$. Then $P_n(X)$ splits modulo~$\ell$ if and only if $c^2-4q$ is a quadratic residue
modulo~$\ell$. In particular, if $P_n(X)$ splits modulo~$\ell$ for \emph{some}~$n\in\NN$, then $P_n(X)$ splits
modulo~$\ell$ for \emph{any}~$n\in\NN$.
\end{lemma}

\begin{proof}
Let $M_1=\left[\begin{smallmatrix} 0 & -q \\ 1 & c\end{smallmatrix}\right]$, and write
    \begin{equation}
    M_1^n=\begin{bmatrix} m_{11} & m_{12} \\ m_{21} & m_{22} \end{bmatrix}.
    \end{equation}
Since $M_1^nM_1=M_1M_1^n$, it follows that $m_{12}=-qm_{21}$ and $m_{22}=m_{11}+cm_{21}$. But then $P_n(X)\equiv
(X-1)(X-q^n)F_n(X)\pmod{\ell}$, where
    \begin{equation}
    F_n(X)\equiv X^2-(2m_{11}+cm_{21})X+m_{11}^2+qm_{21}^2+cm_{11}m_{21}\pmod{\ell}.
    \end{equation}
The discriminant of $F_n(X)$ is given by $\Delta\equiv (c^2-4q)m_{21}^2\pmod{\ell}$; hence the lemma.
\end{proof}

\begin{theorem}\label{theorem:WeilPolynomialSplits}
The Weil polynomial of the Jacobian~$\jac$ of a curve~$C\in\C$ splits modulo~$\ell$.
\end{theorem}

\begin{proof}
For some $n\in\NN$, $\jac[\ell]\subseteq\jac(\FF_{q^n})$. But then~$\frob^n$ acts as the identity on~$\jac[\ell]$, i.e.
$P_n(X)\equiv (X-1)^4\pmod{\ell}$. In particular, $P_n(X)$ splits modulo~$\ell$. But then $P(X)$ splits modulo~$\ell$
by Lemma~\ref{lemma:quadraticResidue}.
\end{proof}

We are now ready to prove the desired result.

\begin{theorem}\label{theorem:FrobeniusDiagonalizable}
The $q$-power Frobenius endomorphism on the Jacobian~$\jac$ of a curve~$C\in\C$ is diagonalizable on~$\jac[\ell]$.
\end{theorem}

\begin{proof}
Cf. Theorem~\ref{theorem:WeilPolynomialSplits}, we may write the Weil polynomial of~$\jac$ as
    \begin{equation}
    P(X)\equiv (X-1)(X-q)(X-\alpha)(X-q/\alpha)\pmod{\ell}.
    \end{equation}
If $\alpha\not\equiv 1,q,q/\alpha\pmod{\ell}$, then the theorem follows. If $\alpha\equiv 1,q\pmod{\ell}$, then
    \begin{equation}
    P(X)\equiv (X-1)^2(X-q)^2\pmod{\ell};
    \end{equation}
in this case, the theorem follows by the last part of the proof of Lemma~\ref{lemma:FrobMatrix}.

Assume that $\alpha\equiv q/\alpha\pmod{\ell}$, i.e. that $\alpha^2\equiv q\pmod{\ell}$. Then the $q$-power Frobenius
endomorphism is represented on~$\jac[\ell]$ by a matrix of the form
    \begin{equation}
    M=\begin{bmatrix}
        1 & 0 & 0 & 0 \\
        0 & q & 0 & 0 \\
        0 & 0 & \alpha & \beta \\
        0 & 0 & 0 & \alpha
      \end{bmatrix}
    \end{equation}
with respect to an appropriate basis of~$\jac[\ell]$. Notice that
    \begin{equation}
    M^{2k}=\begin{bmatrix}
        1 & 0 & 0 & 0 \\
        0 & 1 & 0 & 0 \\
        0 & 0 & 1 & 2k\alpha^{2k-1}\beta \\
        0 & 0 & 0 & 1
      \end{bmatrix}.
    \end{equation}
Thus, $P_{2k}(X)\equiv (X-1)^4\pmod{\ell}$. By Theorem~\ref{theorem:rank:supplement}, it follows that
$\jac[\ell]\subseteq\jac(\FF_{q^{2k}})$. But then $M^{2k}=I$, i.e. $\beta\equiv 0\pmod{\ell}$. Hence, the $q$-power
Frobenius endomorphism on~$\jac$ is diagonalizable on~$\jac[\ell]$ also in this case. The theorem is proved.
\end{proof}

\section{Anti-symmetric pairings on the Jacobian}\label{sec:WeilMatrix}

On $\jac[\ell]$, a non-degenerate, bilinear, anti-symmetric and Galois-invariant pairing
    \begin{equation}\label{eq:Weil-typePairing}
    \pairing:\jac[\ell]\times\jac[\ell]\to\mu_\ell=\grp{\zeta}\subseteq\FF_{q^k}^\times.
    \end{equation}
exists, e.g. the Weil-pairing. Here, $\mu_\ell$ is the group of $\ell^{\textrm{th}}$ roots of unity. Since $\pairing$
is bilinear, it is given by
    \begin{equation}\pairing(x,y)=\zeta^{x^T\pairingMatrix{} y},\end{equation}
for some matrix $\pairingMatrix{}\in\Mat_4(\ZZ/\ell\ZZ)$ with respect to a basis $\basis=\{x_1,x_2,x_3,x_4\}$
of~$\jac[\ell]$. Let~$\frob$ denote the $q$-power Frobenius endomorphism on~$\jac$. Since $\pairing$ is
Galois-invariant,
    \begin{equation}\forall x,y\in\jac[\ell]: \pairing(x,y)^q=\pairing(\frob(x),\frob(y)).\end{equation}
This is equivalent to
    \begin{equation}\forall x,y\in\jac[\ell]: q(x^T\pairingMatrix{}y)=(Mx)^T\pairingMatrix{}(My),\end{equation}
where $M$ is the matrix representation of~$\frob$ on~$\jac[\ell]$ with respect to~$\basis$. Since
$(Mx)^T\pairingMatrix{}(My)=x^TM^T\pairingMatrix{}My$, it follows that
    \begin{equation}\forall x,y\in\jac[\ell]:x^Tq\pairingMatrix{}y=x^TM^T\pairingMatrix{}My,\end{equation}
or equivalently, that~$q\pairingMatrix{}=M^T\pairingMatrix{}M$.

Now, let $\pairing(x_i,x_j)=\zeta^{a_{ij}}$. By anti-symmetry,
    \begin{equation}
    \pairingMatrix{}=
        \begin{bmatrix}
        0        & a_{12}  & a_{13}  & a_{14} \\
        -a_{12} & 0        & a_{23}  & a_{24} \\
        -a_{13} & -a_{23} & 0        & a_{34} \\
        -a_{14} & -a_{24} & -a_{34} & 0
        \end{bmatrix}.
    \end{equation}
Assume that $\frob$ is represented by a diagonal matrix $\diag(1,q,\alpha,q/\alpha)$ with respect to~$\basis$. Then it
follows from $M^T\mathcal{E}M=q\mathcal{E}$, that
    \begin{equation}
    a_{13}(\alpha-q)\equiv a_{14}(\alpha-1)\equiv a_{23}(\alpha-1)\equiv a_{24}(\alpha-q)\equiv 0\pmod{\ell}.
    \end{equation}
If $\alpha\equiv 1,q\pmod{\ell}$, then $\jac(\FF_q)[\ell]$ is bi-cyclic. Hence the following theorem holds.

\begin{theorem}\label{theorem:anti-symmetric-pairings}
Consider a curve $C\in\C$. Let~$\frob$ be the $q$-power Frobenius endomorphism on the Jacobian~$\jac$. Now choose a
basis~$\basis$ of~$\jac[\ell]$, such that~$\frob$ is represented by a diagonal matrix~$\diag(1,q,\alpha,q/\alpha)$ with
respect to~$\basis$. If the $\FF_q$-rational subgroup~$\jac(\FF_q)[\ell]$ of points on the Jacobian of order~$\ell$ is
cyclic, then all non-degenerate, bilinear, anti-symmetric and Galois-invariant pairings on~$\jac[\ell]$ are given by
the matrices
    \begin{equation}
    \pairingMatrix{a,b}=
        \begin{bmatrix}
            0 & a & 0 & 0 \\
            -a & 0 & 0 & 0 \\
            0 & 0 & 0 & b \\
            0 & 0 & -b & 0
        \end{bmatrix},
    \qquad a,b\in(\ZZ/\ell\ZZ)^\times
    \end{equation}
with respect to~$\basis$.
\end{theorem}

\begin{remark}\label{remark:a=b=1}
Let notation and assumptions be as in Theorem~\ref{theorem:anti-symmetric-pairings}. Let~$\pairing$ be a
non-degenerate, bilinear, anti-symmetric and Galois-invariant pairing on~$\jac[\ell]$, and let~$\pairing$ be given
by~$\pairingMatrix{a,b}$ with respect to a basis~$\{x_1,x_2,x_3,x_4\}$ of~$\jac[\ell]$. Then~$\pairing$ is given
by~$\pairingMatrix{1,1}$ with respect to~$\{a^{-1}x_1,x_2,b^{-1}x_3,x_4\}$.
\end{remark}

\begin{remark}\label{remark:AlmostAlwaysCyclic}
In most cases relevant to cryptography, we consider a prime divisor~$\ell$ of size~$q^2$. Assume~$\ell$ is of
size~$q^2$. Then $\ell$ divides neither~$q$ nor~$q-1$. The number of $\FF_q$-rational points on the Jacobian is
approximately~$q^2$. Thus, $\jac(\FF_q)[\ell]$~is cyclic in most cases relevant to cryptography.
\end{remark}

\section{Generators of $\jac[\ell]$}

Consider a curve $C\in\C$ with Jacobian~$\jac$. Assume the $\FF_q$-rational subgroup~$\jac(\FF_q)[\ell]$ of points on
the Jacobian of order~$\ell$ is cyclic. Let~$\frob$ be the $q$-power Frobenius endomorphism on~$\jac$. Let~$\pairing$
be a non-degenerate, bilinear, anti-symmetric and Galois-invariant pairing
    \begin{equation}
    \pairing:\jac[\ell]\times\jac[\ell]\to\mu_\ell=\grp{\zeta}\subseteq\FF_{q^k}^\times.
    \end{equation}
We consider the cases $\ell\nmid 4\tau_k$ and $\ell\mid 4\tau_k$ separately.

\subsection{The case $\ell\nmid 4\tau_k$}

If $\ell$ does not divide~$4\tau_k$, then $\jac(\FF_{q^k})[\ell]$ is bicyclic; cf.~Theo\-rem~\ref{theorem:rank}. Choose
a random point~$\zero\neq x_1\in\jac(\FF_q)[\ell]$, and expand $\{x_1\}$ to a basis $\{x_1,y_2\}$
of~$\jac(\FF_{q^k})[\ell]$, where~$\frob(y_2)=qy_2$. Let $x_2'\in\jac(\FF_{q^k})[\ell]\setminus\jac(\FF_q)[\ell]$ be a
random point. Write $x_2'=\alpha_1x_1+\alpha_2y_2$. Then
    \begin{equation}x_2=x_2'-\frob(x_2')=\alpha_2(1-q)y_2\in\grp{y_2},\end{equation}
i.e. $\frob(x_2)=qx_2$. Now, let~$\jac[\ell]\simeq\jac(\FF_{q^k})[\ell]\oplus W$, where~$W$ is a $\frob$-invariant
submodule of rank two. Choose a random point~$x_3'\in\jac[\ell]\setminus\jac(\FF_{q^k})[\ell]$. Then
    \begin{equation}x_3=x_3'-\frob^k(x_3')\in W\end{equation}
as above. Notice that
    \begin{equation}\jac[\ell]=\grp{x_1,x_2,x_3,\frob(x_3)}\quad\text{if and only if}\quad\pairing(x_3,\frob(x_3))\neq 1;\end{equation}
cf.~Theorem~\ref{theorem:anti-symmetric-pairings}.

Assume~$\pairing(x_3,\frob(x_3))=1$. Then $x_3$ is an eigenvector of~$\frob$. Expand $\{x_1,x_2,x_3\}$ to a basis
$\basis=\{x_1,x_2,x_3,x_4\}$ of~$\jac[\ell]$, such that~$\frob$ is represented by a diagonal matrix on~$\jac[\ell]$
with respect to~$\basis$. We may assume that $\pairing$ is given by~$\pairingMatrix{1,1}$ with respect to~$\basis$;
cf.~Remark~\ref{remark:a=b=1}.

Now, choose a random point~$x\in\jac[\ell]\setminus\jac(\FF_{q^k})[\ell]$. Write
$x=\alpha_1x_1+\alpha_2x_2+\alpha_3x_3+\alpha_4x_4$. Then~$\pairing(x_3,x)=\zeta^{\alpha_4}$. So $\pairing(x_3,x)\neq
1$ if and only if $\ell$ does not divide~$\alpha_4$. On the other hand, $\{x_1,x_2,x_3,x\}$ is a basis of~$\jac[\ell]$
if and only $\ell$ does not divide~$\alpha_4$. Hence, $\{x_1,x_2,x_3,x\}$ is a basis of~$\jac[\ell]$ if and only if
$\ell$ does not divide~$\alpha_4$. Thus, if~$\ell$ does not divide~$4\tau_k$, then the following
Algorithm~\ref{algorithm:NotDivide} outputs generators of~$\jac[\ell]$ with probability~$1-\nicefrac{1}{\ell^n}$.

\begin{algorithm}\label{algorithm:NotDivide}
The following algorithm takes as input a $\C$-curve $C$, the numbers~$\ell$, $q$, $k$ and~$\tau_k$ and a number
$n\in\NN$.
 \begin{enumerate}
  \item Choose points $\zero\neq x_1\in\jac(\FF_q)[\ell]$, $x_2\in\jac(\FF_{q^k})[\ell]\setminus\jac(\FF_q)[\ell]$ and~$x_3'\in
  U:=\jac[\ell]\setminus\jac(\FF_{q^k})[\ell]$; compute~$x_3=x_3'-\frob^k(x_3')$. If $\pairing(x_3,\frob(x_3))\neq 1$,
  then output~$\{x_1,x_2,x_3,\frob(x_3)\}$ and stop.
  \item Let $i=j=0$. While $i<n$ do the following
  \begin{enumerate}
   \item Choose a random point~$x_4\in U$.
   \item $i:=i+1$.
   \item If $\pairing(x_3,x_4)=1$, then $i:=i+1$. Else $i:=n$ and $j:=1$.
  \end{enumerate}
  \item If $j=0$ then output ``failure''. Else output~$\{x_1,x_2,x_3,x_4\}$.
 \end{enumerate}
\end{algorithm}

\subsection{The case $\ell\mid 4\tau_k$}

Assume~$\ell$ divides~$4\tau_k$. Then~$\jac[\ell]\subseteq\jac(\FF_{q^k})$; cf.~Theorem~\ref{theorem:rank:supplement}.
Choose a random point~$\zero\neq x_1\in\jac(\FF_q)[\ell]$, and let $y_2\in\jac[\ell]$ be a point
with~$\frob(y_2)=qy_2$. Write $\jac[\ell]=\grp{x_1,y_2}\oplus W$, where~$W$ is a $\frob$-invariant submodule of rank
two; cf.~the~proof of Lemma~\ref{lemma:FrobMatrix}. Let $\{y_3,y_4\}$ be a basis of~$W$, such that~$\frob$ is
represented on~$\jac[\ell]$ by a diagonal matrix~$M=\diag(1,q,\alpha,q/\alpha)$ on~$\jac[\ell]$ with respect to the
basis
    \begin{equation}
    \basis=\{x_1,y_2,y_3,y_4\}.
    \end{equation}

Now, choose a random point $z\in\jac[\ell]\setminus\jac(\FF_q)[\ell]$. Since $z-\frob(z)\in\grp{y_2,y_3,y_4}$, we may
assume that $z\in\grp{y_2,y_3,y_4}$. Write $z=\alpha_2y_2+\alpha_3y_3+\alpha_4y_4$. Then
    \begin{align*}
    qz-\frob(z)
        &= \alpha_2qy_2+\alpha_3qy_3+\alpha_4qy_4-\left(\alpha_2qy_2+\alpha_3\alpha y_3+\alpha_4(q/\alpha) y_4\right) \\
        &= \alpha_3(q-\alpha)y_3+\alpha_4(q-q/\alpha)y_4;
    \end{align*}
so $qz-\frob(z)\in\grp{y_3,y_4}$. If $qz-\frob(z)=0$, then it follows that $q\equiv 1\pmod{\ell}$. This contradicts the
choice of the curve $C\in\C$. Hence, we have a procedure to choose a point~$\zero\neq w\in W$.

Choose two random points~$w_1,w_2\in W$. Write $w_i=\alpha_{i3}y_3+\alpha_{i4}y_4$ for $i=1,2$. We may assume that
$\pairing$ is given by~$\pairingMatrix{1,1}$ with respect to~$\basis$; cf.~Remark~\ref{remark:a=b=1}. But then
    \begin{equation}\pairing(w_1,w_2)=\zeta^{\alpha_{13}\alpha_{24}-\alpha_{14}\alpha_{23}}.\end{equation}
Hence, $\pairing(w_1,w_2)=1$ if and only if $\alpha_{13}\alpha_{24}\equiv\alpha_{14}\alpha_{23}\pmod{\ell}$.
If~$\alpha_{13}\not\equiv 0\pmod{\ell}$, then $\pairing(w_1,w_2)=1$ if and only if
$\alpha_{24}\equiv\frac{\alpha_{14}\alpha_{23}}{\alpha_{13}}\pmod{\ell}$. So $\pairing(w_1,w_2)\neq 1$ with probability
$1-\nicefrac{1}{\ell}$. Hence, we have a procedure to find a basis of~$W$.

Until now, we have found points $x_1\in\jac(\FF_q)[\ell]$ and $w_3,w_4\in W$, such that $W=\grp{w_3,w_4}$. Now, choose
a random point $x_2\in\jac[\ell]$. Write $x_2=\alpha_1x_1+\alpha_2y_2+\alpha_3y_3+\alpha_4y_4$.
Then~$\pairing(x_1,x_2)=\zeta^{\alpha_2}$, i.e. $\pairing(x_1,x_2)=1$ if and only if $\alpha_2\equiv 0\pmod{\ell}$.
Thus, with probability $1-\nicefrac{\ell^3}{\ell^4}=1-\nicefrac{1}{\ell}$, the set $\{x_1,x_2,w_3,w_4\}$ is a basis
of~$\jac[\ell]$.

Summing up, if~$\ell$ divides~$4\tau_k$, then the following Algorithm~\ref{algorithm:NotDivide} outputs generators
of~$\jac[\ell]$ with probability~$(1-\nicefrac{1}{\ell^n})^2$.

\begin{algorithm}\label{algorithm:Divide}
The following algorithm takes as input a $\C$-curve $C$, the numbers~$\ell$, $q$, $k$ and~$\tau_k$ and a number
$n\in\NN$.
 \begin{enumerate}
  \item Choose a random point $\zero\neq x_1\in\jac(\FF_q)[\ell]$
  \item Let $i=j=0$. While $i<n$ do the following
  \begin{enumerate}
   \item Choose random points $y_3,y_4\in\jac[\ell]$; compute $x_\nu:=q(y_\nu-\frob(y_\nu))-\frob(y_\nu-\frob(y_\nu))$ for $\nu=3,4$.
   \item If $\pairing(x_3,x_4)=1$ then $i:=i+1$. Else $i:=n$ and $j:=1$.
  \end{enumerate}
  \item If $j=0$ then output ``failure'' and stop.
  \item Let $i=j=0$. While $i<n$ do the following
  \begin{enumerate}
   \item Choose a random point $x_2\in\jac[\ell]$.
   \item If $\pairing(x_1,x_2)=1$ then $i:=i+1$. Else $i:=n$ and $j:=1$.
  \end{enumerate}
  \item If $j=0$ then output ``failure''. Else output $\{x_1,x_2,x_3,x_4\}$.
 \end{enumerate}
\end{algorithm}

\subsection{The complete algorithm}

Combining Algorithm~\ref{algorithm:NotDivide} and~\ref{algorithm:Divide} yields the desired algorithm to find
generators of~$\jac[\ell]$.

\begin{algorithm}\label{algorithm:FindGenerators}
The following algorithm takes as input a $\C$-curve $C$, the numbers~$\ell$, $q$, $k$ and~$\tau_k$ and a number
$n\in\NN$.
    \begin{enumerate}
     \item If $\ell\nmid\tau_k$, run Algorithm~\ref{algorithm:NotDivide} on input $(C,\ell,q,k,\tau_k,n)$.
     \item If $\ell\mid\tau_k$, run Algorithm~\ref{algorithm:Divide} on input $(C,\ell,q,k,\tau_k,n)$.
    \end{enumerate}
\end{algorithm}

\begin{theorem}\label{theorem:FindGenerators}
Let $C$ be a $\C$-curve. On input~$(C,\ell,\tau_k,n)$, Algorithm~\ref{algorithm:FindGenerators} outputs generators
of~$\jac[\ell]$ with probability at least~$(1-\nicefrac{1}{\ell^n})^2$ and in expected running time~$O(\log\ell)$.
\end{theorem}

\begin{proof}
We may assume that the time necessary to perform an addition of two points on the Jacobian, to multiply a point with a
number or to evaluate the $q$-power Frobenius endomorphism on the Jacobian is small compared to the time necessary to
compute the (Weil-) pairing of two points on the Jacobian. By~\cite{frey-ruck}, the pairing can be evaluated in
time~$O(\log\ell)$. Hence, the expected running time of Algorithm~\ref{algorithm:FindGenerators} is of
size~$O(\log\ell)$.
\end{proof}

\section{Implementation issues}

A priori, to implement Algorithm~\ref{algorithm:FindGenerators}, we need to find a $q^k$-Weil number~$\omega_k$ of the
Jacobian~$\jac$, in order to check if~$\ell$ ramifies in~$\QQ(\omega_k)$ in the case when $\ell$ divides~$4\tau_k$. On
Jacobians generated by the \emph{complex multiplication method}~\cite{weng03,gaudry,eisen-lauter}, we know the Weil
numbers in advance. Hence, Algorithm~\ref{algorithm:FindGenerators} is particularly well suited for such Jacobians.

Fortunately, in most cases~$\ell$ does not divide~$4\tau_k$, and then we do not have to find a $q^k$-Weil number. And
in fact, we do not even have to compute~$4\tau_k$. To see this, notice that by
Theorem~\ref{theorem:WeilPolynomialSplits}, the Weil polynomial of~$\jac$ is of the form
    \begin{equation}
    P(X)\equiv (X-1)(X-q)(X-\alpha)(X-q/\alpha)\pmod{\ell}.
    \end{equation}
Let~$\frob$ be the $q$-power Frobenius endomorphism on~$\jac$, and let~$P_k(X)$ be the characteristic polynomial
of~$\frob^k$. Since~$\frob$ is diagonalizable on~$\jac[\ell]$, it follows that
    \begin{equation}
    P_k(X)\equiv (X-1)^2(X-\alpha^k)(X-1/\alpha^k)\pmod{\ell}.
    \end{equation}
If $\ell$ divides~$4\tau_k$, then $\jac[\ell]\subseteq\jac(\FF_{q^k})$; cf. Theorem~\ref{theorem:rank:supplement}. But
then $P_k(X)\equiv (X-1)^4\pmod{\ell}$. Hence,
    \begin{equation}\label{check:1}
    \text{$\ell$ divides~$4\tau_k$ if and only if $\alpha^k\equiv 1\pmod{\ell}$.}
    \end{equation}
Assume $\alpha^k\equiv 1\pmod{\ell}$. Then $P_k(X)\equiv (X-1)^4\pmod{\ell}$. Hence,
    \begin{equation}\label{check:2}
    \text{$\ell$ ramifies in~$\QQ(\omega^k)$ if and only if~$\omega^k\notin\ZZ$;}
    \end{equation}
cf.~\cite[Proposition~8.3, p.~47]{neukirch}. Here, $\omega$ is a $q$-Weil number of~$\jac$.

Consider the case when $\alpha^k\equiv 1\pmod{\ell}$ and~$\omega^k\in\ZZ$. Then~$\omega=\sqrt{q}e^{\frac{in\pi}{k}}$
for some $n\in\ZZ$ with $0<n<k$. Assume $k$ divides $mn$ for some $m<k$. Then $\omega^{2m}=q^m\in\ZZ$. Since the
$q$-power Frobenius endomorphism is the identity on the $\FF_q$-rational points on the Jacobian, it follows that
$\omega^{2m}\equiv 1\pmod{\ell}$. Hence, $q^m\equiv 1\pmod{\ell}$, i.e. $k$~divides~$m$. This is a contradiction.
So~$n$ and~$k$ has no common divisors. Let~$\xi=\omega^2/q=e^{\frac{in2\pi}{k}}$. Then~$\xi$ is a primitive
$k^{\text{th}}$ root of unity, and~$\QQ(\xi)\subseteq K$. Since~$[K:\QQ]\leq 4$ and $[\QQ(\xi):\QQ]=\phi(k)$,
where~$\phi$ is the Euler phi function, it follows that~$k\leq 12$. Hence,
    \begin{equation}\label{check:3}
    \text{if $\alpha^k\equiv 1\pmod{\ell}$, then $\omega^k\in\ZZ$ if and only if $k\leq 12$.}
    \end{equation}
The criteria~\eqref{check:1}, \eqref{check:2} and~\eqref{check:3} provides the following efficient
Algorithm~\ref{algorithm:checkCurve} to check whether a given curve is of type $\C$, and whether~$\ell$
divides~$4\tau_k$.

\begin{algorithm}\label{algorithm:checkCurve}
Let~$\jac$ be the Jacobian of a genus two curve~$C$. Assume the odd prime number~$\ell$ divides the number of
$\FF_q$-rational points on~$\jac$, and that~$\ell$ divides neither~$q$ nor~$q-1$. Let~$k$ be the multiplicative order
of~$q$ modulo~$\ell$.
 \begin{enumerate}
  \item Compute the Weil polynomial~$P(X)$ of~$\jac$. Let $P(X)\equiv\prod_{i=1}^4(X-\alpha_i)\pmod{\ell}$.
  \item If $\alpha_i^k\not\equiv 1\pmod{\ell}$ for an $i\in\{1,2,3,4\}$, then output
  ``$C\in\C$ and~$\ell$ does not divide~$4\tau_k$'' and stop.
  \item If $k>12$ then output ``$C\notin\C$'' and stop.
  \item Output ``$C\in\C$ and~$\ell$ divides~$4\tau_k$'' and stop.
 \end{enumerate}
\end{algorithm}

\bibliographystyle{plain}
\bibliography{references}

\end{document}